\documentclass{amsart}

\def\Q{{\mathbf Q}}
\def\Z{{\mathbf Z}}

\def\F{\mathbf{F}}
\def\Gal{\mathrm{Gal}}
\def\End{\mathrm{End}}
\def\Aut{\mathrm{Aut}}
\def\Hom{\mathrm{Hom}}
\def\I{{\mathcal I}}
\def\J{{\mathcal J}}

\def\Sp{\mathrm{Sp}}
\def\M{\mathrm{M}}
\def\dim{\mathrm{dim}}
\def\O{{\mathcal O}}

\def\bmu{\boldsymbol \mu}
\newtheorem{thm}{Theorem}[section]
\newtheorem{lem}[thm]{Lemma}
\newtheorem{cor}[thm]{Corollary}
\newtheorem{prop}[thm]{Proposition}
\theoremstyle{definition}
\newtheorem{defn}[thm]{Definition}

\newtheorem{rem}[thm]{Remark}
\newtheorem{rems}[thm]{Remarks}

\title[Reduction of abelian varieties]
{Reduction of abelian varieties}

\author[A.\ Silverberg]{A.\ Silverberg}
\address{Department of Mathematics, Ohio State University, 
231 W.\ 18 Avenue,
Columbus, Ohio 43210--1174, USA}
\email{silver\char`\@math.ohio-state.edu}

\author[Yu. G. Zarhin]{Yu. G. Zarhin}
\address{Department of Mathematics, Pennsylvania State University, 
University Park, PA 16802, USA}
\email{zarhin\char`\@math.psu.edu}

\begin{document}

\maketitle

\section{Introduction}

In this paper we study the reduction of abelian varieties.
We assume $F$ is a field with a discrete valuation $v$, $X$ is an
abelian variety over $F$, and $n$ is an integer not divisible by
the residue characteristic.

In Part \ref{semistabpart} we give criteria for semistable reduction.
Suppose $n \ge 5$. In Theorem \ref{ssredlem} 
we show that $X$ has semistable reduction if and only if 
$(\sigma - 1)^2 = 0$ on the $n$-torsion in $X$, for every $\sigma$ in 
the absolute inertia group.
In Theorem \ref{ssredconverse} we show (using Theorem \ref{ssredlem}) 
that $X$ has semistable reduction 
if and only if there exists a subgroup of $n$-torsion points such that
the absolute inertia group acts trivially on both it and its orthogonal
complement with respect to the $e_n$-pairing.
We deduce as special cases both Raynaud's criterion (that the abelian 
variety have full level $n$ structure for $n \ge 3$; see 
Theorem \ref{raynaud}) and the criterion of \cite{semistab} (that the 
abelian variety have partial level $n$ structure for $n \ge 5$; see 
Theorem \ref{ssred}). We also obtain a (near) converse to the
criterion of \cite{semistab}. The proofs are based on the fundamental 
results of Grothendieck on semistable reduction of abelian varieties 
(see \cite{SGA}). In \S\ref{except} we allow $n<5$. In
\S\ref{Gsect} we give a measure of potentially good reduction.
We discuss other measures of potentially good reduction in 
Part  \ref{neronpart}.

In Part \ref{neronpart} we study N\'eron models of abelian 
varieties with potentially good reduction and torsion points of 
small order. Suppose that the valuation ring is henselian and 
the residue field is algebraically closed. 
If $X$ has good reduction, 
then $X_n \subseteq X(F)$ (this is an immediate corollary of the existence of 
N\'eron models; see Lemma \ref{neronlemma} below). 
On the other hand, if $X_n \subseteq X(F)$ and $n \ge 3$, 
then by virtue of Raynaud's criterion for semistable reduction, 
$X$ has good reduction.
Notice that the failure of $X$ to have good reduction is measured by 
the dimension $u$ of the unipotent radical of the special fiber of the N\'eron 
model of $X$. In particular, $u = 0$ if and only if $X$ has good reduction. 
In general, $0 \le u \le \dim(X)$. The equality $u = \dim(X)$
says that $X$ has purely additive reduction.
Another measure of the deviation from good reduction is the (finite) group of
connected components $\Phi$ of the special fiber of the N\'eron model. 
If $X$ has good reduction then $\Phi=\{0\}$, but the converse statement
is not true in general.

The aim of \S\ref{main} is to connect explicitly the invariants 
$u$ and $\Phi$ with the failure of $X(F)$ to contain all the 
$n$-torsion points. This failure can be
measured by the index $[X_n:X_n(F)]$. We assume that at 
least ``half'' of the
$n$-torsion points are rational over $F$. More precisely, we assume that 
there exists an $F$-rational polarization $\lambda$ on $X$ 
and a maximal isotropic 
(with respect to the pairing $e_{\lambda,n}$ induced from the 
Weil $e_n$-pairing by $\lambda$) 
subgroup of $X_n$ consisting of $F$-rational points.
If in addition $n \ge 5$, then $X$ has good reduction (see Theorem 7.4 
of \cite{semistab}), and therefore $u = 0$, $\Phi=\{0\}$, and $X_n=X_n(F)$. 
Therefore, we have to investigate only the cases $n = 2$, $3$, and $4$.
Let $\Phi'$ denote the prime-to-$p$ part of $\Phi$, where $p$ is the
residue characteristic (with $\Phi ' = \Phi$ if $p = 0$). 

We show that if $n = 2$ then $\Phi'$ is an elementary abelian
$2$-group and  $[X_2:X_2(F)]\#\Phi' = 4^{u}$, 
if $n = 3$ then $[X_3:X_3(F)] = 3^{u}$ and
$\Phi' \cong (\Z/3\Z)^u$, and
if $n = 4$ then $X_2 \subseteq X(F)$,
$[X_4:X_4(F)] = 4^{u}$, and $\Phi' \cong (\Z/2\Z)^{2u}$.
If instead of assuming partial level $n$ structure we assume that 
all the points of order $2$ on $X$ are defined over $F$, then
$[X_4:X_4(F)]=4^u$ and $\Phi' \cong (\Z/2\Z)^{2u}$. 

Earlier work on abelian varieties with potentially good reduction and
on groups of connected components of N\'eron models 
has been done by Serre and Tate \cite{Serre-Tate},
Silverman \cite{Silverman},
Lenstra and Oort \cite{LenstraOort}, 
Lorenzini \cite{Lorenzini}, and Edixhoven \cite{Edixhoven}. 

Silverberg would like to thank the IHES and the Bunting
Institute for their hospitality, and the NSA and the Science Scholars
Fellowship Program at the Bunting Institute for financial support.
Zarhin would like to thank the NSF for financial support. 
He also would like to thank the organizers of the 
NATO/CRM 1998 Summer School 
on the Arithmetic and Geometry of Algebraic Cycles for 
twelve wonderful days in Banff.

\section{Notation and definitions}

If $F$ is a field, let $F^s$ denote a separable closure. 
Suppose that $X$ is an abelian variety defined over 
$F$, and $n$ is a positive integer not divisible by the characteristic
of $F$. Let $X^\ast$ denote the dual abelian variety of $X$, let
$X_n$ denote the kernel of multiplication by $n$ in $X(F^s)$, let
$X_n^\ast$ denote the kernel of multiplication by $n$ in $X^\ast(F^s)$,
and let $\bmu _n$ denote the $\Gal (F^s/F)$-module
of $n$-th roots of unity in $F^s$. 
The $e_n$-pairing 
$$e_n : X_n \times X_n^\ast \to {\bmu}_n$$
is a $\Gal(F^s/F)$-equivariant nondegenerate pairing 
(see \S 74 of \cite{WeilAV}).
If $S$ is a subgroup of $X_n$, let 
$$S^{\perp_n} = 
\{ y \in X_n^\ast : e_n(x,y) = 1 \text{ for every } x \in S \}
\subseteq X_n^\ast.$$
If $\lambda$ is a polarization on $X$, define
$$e_{\lambda ,n} : X_n \times X_n \to {\bmu}_n$$
by $e_{\lambda ,n}(x,y) = e_n(x,\lambda(y))$
(see \S 75 of \cite{WeilAV}). Then
$$\sigma (e_{\lambda ,n}(x_1,x_2)) = 
e_{\sigma (\lambda ),n}(\sigma (x_1),\sigma (x_2))$$ for every 
$\sigma \in \Gal (F^s/F)$ and $x_1$, $x_2 \in X_n$. If $n$ is relatively prime 
to the degree of the polarization $\lambda$, then the pairing 
$e_{\lambda ,n}$ is nondegenerate.
If $\ell$ is a prime not equal to the characteristic of $F$,
and $d = \dim(X)$, let
$$\rho_{\ell,X} : \Gal(F^s/F) \to \Aut(T_\ell(X)) \cong \M_{2d}(\Z_\ell)$$
denote the $\ell$-adic representation on the Tate module $T_\ell(X)$
of $X$, and let $V_\ell(X) = T_\ell(X) \otimes_{\Z_\ell} \Q_\ell$.
Let $I$ denote the identity matrix in $\M_{2d}(\Z_\ell)$.

If $L$ is a Galois extension of $F$,
$v$ is a discrete valuation on $F$, and $w$ is an extension of $v$ to $L$,
let $\I(w/v)$ denote the inertia subgroup at $w$ of $\Gal(L/F)$.
If $X$ is an abelian variety over $F$, let $X_v$ denote the 
special fiber of the N\'eron model
of $X$ at $v$ and let $X_v^0$ denote its identity connected
component. Let $a$, $u$, and $t$ denote, respectively, the abelian,
unipotent, and toric ranks of $X_v$. Then $a + u + t = \dim(X)$.
If $p$ ( $\ge 0$) is the residue characteristic of $v$, 
let $\Phi '$ denote the prime-to-$p$ part of the group of 
connected components of the special
fiber of the N\'eron model of $X$ at $v$ (with $\Phi '$ the full
group of components if $p = 0$). 

\begin{defn}
If $v$ is a discrete valuation on a field $F$, we say
the valuation ring is 
{\em strictly henselian} if the valuation ring is henselian and 
the residue field is algebraically closed.
\end{defn}

\begin{defn}
Suppose $L/F$ is an extension of fields, $w$ is a discrete valuation
on $L$, and $v$ is the restriction of $w$ to $F$. We say that
$w/v$ is {\em unramified} if a uniformizing element of the
valuation ring for $v$
induces a uniformizing element of the valuation ring for $w$
and the residue field extension is separable (see Definition 1
on page 78 of \cite{BLR}).
\end{defn}

\begin{rem}[Remark 5.3 of \cite{semistab}]
\label{ramifiedcyclic}
Suppose $v$ is a discrete valuation on a field $F$, and $m$ is a positive
integer not divisible by the residue characteristic. Then every 
degree $m$ Galois extension of $F$ totally
ramified at $v$ is cyclic. If $F(\zeta_m) = F$,
then $F$ has a cyclic extension of degree $m$ which is totally ramified at
$v$. In particular, if the residue characteristic is not $2$ then 
$F$ has a quadratic extension which is (totally and tamely) ramified at $v$.
If the valuation ring is henselian and the residue field is separably closed, 
then $F= F(\zeta_m)$ and therefore $F$ has a cyclic totally ramified extension 
of degree $m$. (See Remark 5.3 of \cite{semistab}.)
Note also that $F$ has no non-trivial unramified extensions if and only if 
the valuation ring is henselian and the residue field is separably closed.
\end{rem}

\part{Semistable reduction of abelian varieties}
\label{semistabpart}

\section{Preliminaries}

\begin{defn}
If $k$ is a positive integer, define a finite set 
of prime powers $N(k)$ by 
$$N(k) = \{\text{prime powers $\ell^m : 0 \le m(\ell - 1) \le k $}\}.$$
\end{defn}
For example, $$N(1) = \{1, 2\}, \quad N(2) = \{1, 2, 3, 4\},$$
$$N(3) = \{1, 2, 3, 4, 8\},  \quad  
N(4) = \{1, 2, 3, 4, 5, 8, 9, 16\}.$$

\begin{thm}
\label{quasithm}
Suppose $n$ and $k$ are positive integers, $\O$ is an integral 
domain of characteristic zero such that no rational prime which 
divides $n$ 
is a unit in $\O$, $\alpha \in \O$, $\alpha$ has finite 
multiplicative order, 
and $(\alpha - 1)^k \in n\O$.
If $n \notin N(k)$, then $\alpha = 1$. In particular, if
$(\alpha - 1)^2 \in n\O$ and $n \ge 5$, then $\alpha = 1$.
\end{thm}

\begin{proof}
See Corollary 3.3 of \cite{serrelem}.
\end{proof}

\begin{lem}[Lemma 5.2 of \cite{semistab}]
\label{localglobal}
Suppose that $d$ and $n$ are positive integers, and for each prime $\ell$ 
which divides
$n$ we have a matrix $A_\ell \in M_{2d}(\Z_\ell)$ such that the
characteristic polynomials of the $A_\ell$ have integral coefficients
independent of $\ell$, and such that $(A_\ell - I)^2 \in nM_{2d}(\Z_\ell)$.
Then for every eigenvalue $\alpha$ of $A_\ell$, $(\alpha - 1)/\sqrt{n}$
satisfies a monic polynomial with integer coefficients.
\end{lem}

\begin{lem}[Lemma 4.2 of \cite{degree}]
\label{mevals}
Suppose $v$ is a discrete valuation on a field $F$ with residue
characteristic $p \ge 0$, $m$ is a positive integer, $\ell$ is a prime,
$p$ does not divide $m\ell$, $K$ is a degree $m$ extension of $F$
which is totally ramified above $v$, and ${\bar v}$ is an extension of $v$ to
a separable closure $K^s$ of $K$. Suppose that $X$ is an abelian variety over $F$,
and for every $\sigma \in \I({\bar v}/v)$,
all the eigenvalues of $\rho_{\ell,X}(\sigma)$ 
are $m$-th roots of unity. Then $X$ has
semistable reduction at the extension of $v$ to $K$.
\end{lem}

\section{Criteria for semistable reduction}

\begin{thm}[Galois Criterion for Semistable Reduction]
\label{galcrit}
Suppose $X$ is an abel\-ian variety over a field $F$, $v$ is a discrete 
valuation on $F$, $\ell$ is a prime not equal to the residue
characteristic of $v$, ${\bar v}$ is an extension of $v$ to $F^s$,
and $\I = \I({\bar v}/v)$. 
Then the following are equivalent:
\begin{enumerate}
\item[(i)] $X$ has semistable reduction at $v$,
\item[(ii)] $\I$ acts unipotently on $T_\ell(X)$; i.e.,
all the eigenvalues of $\rho_{\ell,X}(\sigma)$ are $1$, 
for every $\sigma \in \I$,
\item[(iii)] for every $\sigma \in \I$, 
$(\rho_{\ell,X}(\sigma) - I)^2 = 0$.
\end{enumerate}
\end{thm}

\begin{proof}
See Proposition 3.5 and Corollaire 3.8 of \cite{SGA} 
and Theorem 6 on p.~184 of \cite{BLR}.
\end{proof}

\begin{thm}[Raynaud Criterion for Semistable Reduction]
\label{raynaud}
\hfil Suppose $X$ is an abelian variety over a field $F$ with a 
discrete valuation $v$, $m$ is a positive integer not divisible by the
residue characteristic of $v$, and the points of $X_m$ are defined
over an extension of $F$ which is unramified over $v$. If $m \ge 3$, then $X$ has
semistable reduction at $v$. 
\end{thm}

\begin{proof}
See Proposition 4.7 of \cite{SGA}.
\end{proof}

\begin{prop}
\label{sslem}
Suppose $X$ is an abelian variety over a field $F$, $v$ is a discrete 
valuation on $F$, $n$ is an integer not 
divisible by the residue characteristic of $v$,  
${\bar v}$ is an extension of $v$ to $F^s$,
and $\I = \I({\bar v}/v)$. Let $S = X_n^\I$, the elements of
$X_n$ on which $\I$ acts as the identity. If
$X$ has semistable reduction at $v$, then
\begin{enumerate}
\item[{(i)}] $(\sigma - 1)^2 X_n = 0$ for every $\sigma \in \I$, and
\item[{(ii)}] $\I$ acts as the identity on $S^{\perp_n}$.
\end{enumerate}
\end{prop}

\begin{proof}
Suppose $X$ has semistable reduction at $v$.
By Theorem \ref{galcrit}, we have (i).
It follows that $\sigma^n = 1$ on $X_n$. Since $n$ is 
not divisible by the residue characteristic, $X_n$ is tamely ramified
over $F$. Let $\J$ denote the first ramification group. Then the
action of $\I$ on $X_n$ factors through $\I/\J$. Let $\tau$ denote
a lift to $\I$ of a topological generator of the pro-cyclic group
$\I/\J$. Since 
$$e_n((\tau - 1)X_n,(X_n^\ast)^\I) = 1,$$
we have 
$$\#((X_n^\ast)^\I)\#((\tau - 1)X_n) \le \#X_n^\ast.$$
The map from $X_n$ to $(\tau - 1)X_n$ defined by
$y \mapsto (\tau - 1)y$ defines a short exact sequence
$$0 \to S \to X_n \to (\tau - 1)X_n \to 0.$$
Therefore, 
$$\#S\#((\tau - 1)X_n) = \#X_n = \#S\#S^{\perp_n}.$$
Similarly, 
$$\#((X_n^\ast)^\I)\#((\tau - 1)X_n^\ast) = \#X_n^\ast.$$
Therefore, 
$$\#S^{\perp_n} = \#((\tau - 1)X_n) \le \#((\tau - 1)X_n^\ast).$$
Since $(\tau - 1)X_n^\ast \subseteq S^{\perp_n}$,
we conclude that 
$$S^{\perp_n} = (\tau - 1)X_n^\ast.$$
By (i),
we have $(\tau - 1)^2X_n = 0$. It follows from the natural
$\Gal(F^s/F)$-equivariant isomorphism 
$X_n^\ast \cong \Hom(X_n,\bmu_n)$ that $(\tau - 1)^2X_n^\ast = 0$,
and therefore $\I$ acts as the identity on $S^{\perp_n}$.
\end{proof}

\begin{thm}
\label{ssredlem}
Suppose $X$ is an abelian variety over a field $F$, $v$ is a discrete 
valuation on $F$, $n$ is an integer not 
divisible by the residue characteristic of $v$, $n \ge 5$, 
${\bar v}$ is an extension of $v$ to $F^s$,
and $\I = \I({\bar v}/v)$. 
Then $X$ has semistable reduction at $v$ if and only if
$(\sigma - 1)^2 X_n = 0$ for every $\sigma \in \I$.
\end{thm}

\begin{proof}
If $X$ has semistable reduction at $v$ then
for every $\sigma \in \I$ we have $(\sigma - 1)^2 X_n = 0$, 
by Proposition \ref{sslem}i. 

Conversely, suppose $n \ge 5$ and $(\sigma - 1)^2 X_n = 0$
for every $\sigma \in \I$.
Let $\I' \subseteq \I$ be the inertia group for the prime 
below ${\bar v}$ in a
finite Galois extension of $F$ over which $X$ has 
semistable reduction.
Take $\sigma \in \I$. Then 
$\sigma^m \in \I'$ for some $m$.
Let $\ell$ be a prime divisor of $n$. 
Theorem \ref{galcrit} implies that
$(\rho_{\ell,X}(\sigma)^m - I)^2 = 0$.
Let $\alpha$ be an eigenvalue of $\rho_{\ell,X}(\sigma)$.
Then $(\alpha^m - 1)^2 = 0$.
Therefore, $\alpha^m = 1$.
By our hypothesis, 
$$(\rho_{\ell,X}(\sigma) - I)^2 \in n\M_{2d}(\Z_\ell),$$
where $d = \dim(X)$.
By Theorem 4.3 on p.~359 of \cite{SGA}, 
the characteristic polynomial of $\rho_{\ell,X}(\sigma)$
has integer coefficients which are independent of $\ell$.
By Lemma \ref{localglobal},
$(\alpha - 1)^2 \in n{\bar \Z}$, where ${\bar \Z}$ denotes
the ring of algebraic integers.
Since $n \ge 5$,
by Theorem \ref{quasithm} we have $\alpha = 1$ (i.e., 
$\I$ acts unipotently on $T_\ell(X)$).
By Theorem \ref{galcrit}, $X$ has semistable reduction at $v$.
\end{proof}

\begin{thm}
\label{ssredconverse}
Suppose $X$ is an abelian variety over a field $F$, $v$ is a discrete 
valuation on $F$, $n$ is an integer not divisible by the residue 
characteristic of $v$, $n \ge 5$, ${\bar v}$ is an extension of
$v$ to $F^s$, and $\I = \I({\bar v}/v)$. Then $X$ has semistable 
reduction at $v$ if and only if there exists a subgroup $S$ of $X_n$
such that $\I$ acts as the identity on $S$ and on $S^{\perp_n}$.
\end{thm}

\begin{proof}
Suppose there exists a subgroup $S$ as in the statement of the
theorem. The map $x \mapsto (y \mapsto e_n(x,y))$ induces a
$\Gal(F^s/F)$-equivariant isomorphism from $X_n/S$ onto
$\Hom(S^{\perp_n},\bmu_n)$. Suppose $\sigma \in \I$. Then 
$\sigma = 1$ on $S^{\perp_n}$ and on $\bmu_n$. Therefore, $\sigma = 1$ on
$X_n/S$. Thus, 
$(\sigma - 1)^2X_n \subseteq (\sigma - 1)S = 0$.
By Theorem \ref{ssredlem}, $X$ has semistable reduction at $v$.

Conversely, suppose $X$ has semistable reduction at $v$. 
Let $S = X_n^\I$, and apply Proposition \ref{sslem}ii.
\end{proof}

\begin{thm}
\label{ssred}
Suppose $X$ is an abelian variety over a field $F$, $v$ is a discrete 
valuation on $F$, $\lambda$ is a polarization on $X$ defined over
an extension of $F$ which is unramified over $v$, 
$n$ is a positive integer not 
divisible by the residue characteristic of $v$, and $n \ge 5$.
\begin{enumerate} 
\item[{(i)}]
If $\widetilde{X}_n$ is a maximal isotropic subgroup of $X_n$ with 
respect to $e_{\lambda,n}$, and the points of 
$\widetilde{X}_n$ are defined over an extension of $F$ which is 
unramified over $v$, then $X$ has semistable reduction at $v$.
\item[{(ii)}] Conversely, if $X$ has semistable reduction at $v$,
and the degree of the polarization $\lambda$ is relatively prime
to $n$, then there exists a maximal isotropic subgroup
of $X_n$ with respect to $e_{\lambda,n}$, all of
whose points are defined over an extension 
of $F$ which is unramified over $v$.
\end{enumerate}
\end{thm}

\begin{proof}
Under the hypotheses in (i),
let $S = \widetilde{X}_n$. Then $S^{\perp_n} = \lambda(S)$, and
$X$ has semistable reduction at $v$ by applying 
Theorem \ref{ssredconverse}.

Conversely, suppose $X$ has semistable reduction at $v$. Let
${\bar v}$ be an extension of $v$ to $F^s$ and let 
$\I = \I({\bar v}/v)$. Let $S = X_n^\I$. If $G$ is a subgroup of
$X_n$, let
$$G^{\perp_{\lambda,n}} = 
\{ y \in X_n : e_{\lambda,n}(x,y) = 1 \text{ for every } x \in G \}.$$
Since the degree of $\lambda$ is relatively prime to $n$, 
$\lambda$ induces an isomorphism between $S^{\perp_{\lambda,n}}$
and $S^{\perp_n}$. Since $\lambda$ is defined over an unramified
extension, $\I$ acts as the identity on $S^{\perp_{\lambda,n}}$
by Proposition \ref{sslem}ii. Therefore, 
$S^{\perp_{\lambda,n}} \subseteq S = X_n^\I$. 
The pairing $e_{\lambda,n}$ induces a nondegenerate pairing on
$S/S^{\perp_{\lambda,n}}$. Let $H$ be the inverse image
in $S$ (under the natural projection) of a maximal isotropic 
subgroup of $S/S^{\perp_{\lambda,n}}$. It is easy to check that 
$H$ is a maximal isotropic subgroup of $X_n$ with respect
to $e_{\lambda,n}$, proving (ii).  
\end{proof}

\begin{rems}
Raynaud's criterion (Theorem \ref{raynaud}) follows from 
Theorem \ref{ssredconverse} by letting 
$n = m^2$ and $S = X_m \subset X_{n}$
(since then $S^{\perp_n} = X^\ast_m$, the dual Galois module
of $X_m$, and $n \ge 5$ whenever $m \ge 3$). 
The converse of Raynaud's criterion is clearly false, 
i.e., semistable reduction does not imply that
the $n$-torsion points are unramified (for $n \ge 3$ and $n$ not
divisible by the residue characteristic), as can be seen, for
example, by
comparing Raynaud's criterion with the N\'eron-Ogg-Shafarevich
criterion for good reduction, and considering an abelian variety
with semistable but not good reduction.

Theorem \ref{ssred}i is Theorem 6.2 of \cite{semistab}.
Similarly, the other results of \cite{semistab} and of \S 3 of
\cite{connected} can readily be generalized to the setting 
of Theorem \ref{ssredconverse}.
Theorem \ref{ssred}ii shows that the sufficient condition for
semistability given in Theorem 6.2 of \cite{semistab} comes close
to being a necessary condition. 
Note that Theorem \ref{ssred}ii would be false if the condition on 
the degree of the polarization were omitted.
\end{rems}

\begin{defn}
Suppose $v$ is a discrete valuation on $F$ of residue characteristic
$p$. 
We say $v$ satisfies (*) if at least one of the following
conditions is satisfied:
\begin{enumerate}
\item[(a)]  $p \ne 2$, 
\item[(b)] the valuation ring 
is henselian and the residue field is separably closed.
\end{enumerate}
\end{defn}

The techniques of the above proofs can be extended to prove
the following result. The proof will appear in \cite{etale}.

\begin{thm}
\label{highercoh}
Suppose $X$ is an abelian variety over a field $F$, and $v$ is a discrete 
valuation on $F$ of residue characteristic $p \ge 0$.
Suppose $k \in \Z$, and
$0 < k < 2\dim(X)$.
\begin{enumerate}
\item[{(i)}] If either $X$ has semistable reduction at $v$, 
or $k$ is even and 
$X$ has purely additive reduction at $v$ which becomes
semistable over a quadratic extension of $F$,
then 
$$(\sigma - 1)^{k+1}H^k_{\text{\'et}}(X \times_F F^s, \Z_{\ell}) = 0$$ 
for every $\sigma \in \I$ and every prime $\ell \ne p$,
and
$$(\sigma - 1)^{k+1}H^k_{\text{\'et}}(X \times_F F^s, \Z/n\Z) = 0$$ 
for every $\sigma \in \I$ and every positive integer $n$ not
divisible by $p$.
\item[{(ii)}] Suppose $n$ is a positive integer not divisible by
$p$, and 
$$(\sigma - 1)^{k+1}H^k_{\text{\'et}}(X \times_F F^s, \Z/n\Z) = 0$$ 
for every $\sigma \in \I$. 
Suppose $L$ is a degree $R(k+1,n)$ extension of $F$ which is
totally ramified above $v$, and let $w$ be the extension of
$v$ to $L$.
If $k$ is odd, then $X$ has semistable reduction at $w$.
If $k$ is even and $v$ satisfies (*), 
then either $X$ has semistable reduction at $w$, or
$X$ has purely additive reduction at $w$ which becomes
semistable over a quadratic extension of $L$.
\end{enumerate}
\end{thm}

If we restrict to the case where $n \notin N(k+1)$, we obtain
the following result. This result gives necessary and sufficient
conditions for semistable reduction, and also necessary and
sufficient conditions for $X$ to have either semistable reduction
or purely additive reduction which becomes semistable over a
quadratic extension.

\begin{cor}
\label{highercohcor}
Suppose $X$ is an abelian variety over a field $F$, $v$ is a discrete 
valuation on $F$ of residue characteristic $p \ge 0$, $k$ and $n$
are positive integers, $\ell$ is a prime number,
$k < 2\dim(X)$, $n$ and $\ell$ are not divisible by $p$, 
and $n \notin N(k+1)$.
\begin{enumerate}
\item[{(i)}] Suppose $k$ is odd. Then the following are
equivalent:
\begin{enumerate}
\item[{(a)}]  $X$ has semistable reduction at $v$,
\item[{(b)}] for every $\sigma \in \I$,
$$(\sigma - 1)^{k+1}H^k_{\text{\'et}}(X \times_F F^s, 
\Z_{\ell}) = 0,$$ 
\item[{(c)}]  for every $\sigma \in \I$,
$$(\sigma - 1)^{k+1}H^k_{\text{\'et}}(X \times_F F^s, \Z/n\Z) = 0.$$ 
\end{enumerate}
\item[{(ii)}] Suppose $k$ is even and $v$ satisfies (*).
Then the following are equivalent:
\begin{enumerate}
\item[{(a)}]  either $X$ has semistable reduction at $v$, or
$X$ has purely additive reduction at $v$ which becomes
semistable over a quadratic extension of $F$,
\item[{(b)}] for every $\sigma \in \I$,
$$(\sigma - 1)^{k+1}H^k_{\text{\'et}}(X \times_F F^s, \Z_{\ell}) = 
0,$$
\item[{(c)}] for every $\sigma \in \I$,
$$(\sigma - 1)^{k+1}H^k_{\text{\'et}}(X \times_F F^s, \Z/n\Z) = 0.$$ 
\end{enumerate}
\end{enumerate}
\end{cor}

\section{Exceptional $n$}
\label{except}
In this section we discuss briefly the ``exceptional'' cases 
$n=2,3,4$. For the proofs, and for examples which show the results
are sharp, we refer the reader to \cite{degree}.

First, let us state the following ``one-way"
generalization of Theorem 4.5.

\begin{thm}
\label{oneway}
Suppose $X$ is an abelian variety over a field $F$, $v$ is a discrete 
valuation on $F$, and $n$ is an integer greater than $1$ which is not 
divisible by the residue characteristic of $v$. 
Suppose there exists a subgroup $S$ of $X_n$
such that $\I$ acts as the identity on $S$ and on $S^{\perp_n}$.
Then $X$ has semistable reduction over every degree $R(n)$
extension of $F$ totally ramified above $v$.
\end{thm}

It turns out that the converse statement is not true. 
However, the following result gives an ``approximate converse''.

\begin{thm}
\label{bothways}
Suppose $n = 2$, $3$, or $4$, respectively. Suppose  
$X$ is an abelian variety over a field $F$, and 
$v$ is a discrete valuation on $F$ whose residue characteristic 
does not divide $n$. 
Suppose
$L$ is an extension of $F$ of degree $4$, $3$, or
$2$, respectively, which 
is totally ramified above $v$. 
Then the following are equivalent:
\begin{enumerate}
\item[(i)] $X$ has semistable reduction over $L$ above $v$,
\item[(ii)] there exist an abelian variety $Y$ over a finite
extension $K$ of $F$ unramified above $v$, a separable
$K$-isogeny $\varphi : X \to Y$, 
and a subgroup 
$S$ of $Y_n$ such that 
$\I$ acts as the identity on $S$ and 
on $S^{\perp_n}$.
\end{enumerate}
Further, $\varphi$ can be taken so that
its kernel is killed by $8$, $9$, or $4$, respectively.
If $X$ has potentially good reduction at $v$, then $\varphi$ 
can be taken so that its kernel is killed by $2$, $3$, or
$2$, respectively.
\end{thm}

In the case of low-dimensional $X$ this result may be improved as follows.

\begin{thm}
\label{ellcor}
In Theorem \ref{bothways},
with $d = \dim(X)$, 
$\varphi$ can be taken so that its kernel is killed
by $4$ if $d = 3$ and $n = 2$,
by $3$ if $d =2$ and $n = 3$, and 
by $2$ if $d = n = 2$. If $d = 1$, then we can take
$Y = X$ and $\varphi$ the identity map.
\end{thm}

In the case of elliptic curves this implies the following statement.

\begin{cor}
\label{4326cor}
Suppose 
$X$ is an elliptic curve over a field $F$, and 
$v$ is a discrete valuation on $F$ of residue characteristic 
$p \ge 0$. 
\begin{enumerate}
\item[(a)] If $p \ne 2$, 
then $X$ has semistable reduction above $v$ over a 
totally ramified quartic extension of $F$  
if and only if 
$X$ has an $\I$-invariant point of order $2$.
\item[(b)] If $p \ne 3$, 
then $X$ has semistable reduction above $v$ over a totally ramified  
cubic extension of $F$  
if and only if $X$ has an $\I$-invariant point of order $3$.
\item[(c)] 
If $p \ne 2$,
then $X$ has semistable reduction above $v$ over a quadratic 
extension of $F$ if and only if 
either $X$ has an $\I$-invariant point of order $4$,
or all the points of order $2$ on $X$ are $\I$-invariant.
\item[(d)] 
If $p \ne 2$ and $X$ has bad but potentially good reduction
at $v$,
then $X$ has good reduction above $v$ over a quadratic 
extension of $F$ if and only if 
$X$ has no $\I$-invariant point of order $4$ 
and all its points of order $2$ are $\I$-invariant.
\item[(e)] Suppose $p$ is not $2$ or $3$. Then the following
are equivalent:
\begin{enumerate} 
\item[(i)]  $X$ has no $\I$-invariant points of order $2$ or $3$, 
\item[(ii)] there does not exist a finite separable extension $L$
of $F$ of degree less than $6$ 
such that $X$ has semistable reduction at the restriction of 
${\bar v}$ to $L$.
\end{enumerate}
\item[(f)] Suppose $p$ is not $2$ or $3$. Then the following
are equivalent:
\begin{enumerate} 
\item[(i)]  $X$ has no $\I$-invariant points of order $4$ or $3$
and not all the points of order $2$ are $\I$-invariant, 
\item[(ii)] there does not exist a finite separable extension $L$
of $F$ of degree less than $4$ 
such that $X$ has semistable reduction at the restriction 
of ${\bar v}$ to $L$.
\end{enumerate}
\end{enumerate}
\end{cor}

	In the case of potentially good reduction the following statement
holds true.

\begin{thm}
\label{paddcor}
Suppose 
$X$ is an abelian variety over a field $F$,  
$v$ is a discrete valuation on $F$ of residue characteristic 
$p \ge 0$, and $X$ has purely additive and  
potentially good reduction at $v$. 
\begin{enumerate}
\item[(a)] 
If $p \ne 2$, then $X$ has good reduction above $v$ over a quadratic 
extension of $F$ if and only if 
there exists a subgroup $S$ of $X_4$ such 
that $\I$ acts as the identity on $S$ and on $S^{\perp_4}$.
\item[(b)] If $p \ne 3$, 
then $X$ has good reduction above $v$ over a totally ramified cubic
extension of $F$ if and only if 
there exists a subgroup $S$ of $X_3$ such 
that $\I$ acts as the identity on $S$ and on $S^{\perp_3}$.
\item[(c)] Suppose $p \ne 2$, and
$L/F$ is a degree $4$ extension, 
totally ramified above $v$, 
which has a quadratic subextension over which $X$ has purely
additive reduction.  
Then $X$ has good reduction above $v$ over $L$ if and only if 
there exists a subgroup $S$ of $X_2$ such 
that $\I$ acts as the identity on $S$ and on $S^{\perp_2}$.
\end{enumerate}
\end{thm}

\section{A measure of potentially good reduction}
\label{Gsect}

Suppose $v$ is a discrete valuation on a field $F$, and
$X$ is an abelian variety over $F$ which has potentially good 
reduction at $v$. Let $F_{v}^{nr}$ denote the maximal unramified
extension of the completion of $F$ at $v$, let $L$ denote
the smallest extension of $F_{v}^{nr}$ over which $X$ has good
reduction, and let
$$G_{v,X} = \Gal(L/F_{v}^{nr}).$$
Then $G_{v,X}$ can also be characterized as the inertia group of
the extension $F(X_n)/F$, where $n$ is any integer 
greater than $2$ and not divisible by the residue characteristic
of $v$ (see Corollary 2 on p.~497 of \cite{Serre-Tate}).
Clearly, $X$ has good reduction at $v$ if and only if
$G_{v,X} = 1$. The finite group $G_{v,X}$ is a measure of how
far $X$ is from having good reduction at $v$.

If $A$ is a matrix, let $P_A$ denote its characteristic polynomial.
The following result gives constraints on the group $G_{v,X}$.

\begin{thm}
\label{Gthm}
Suppose $v$ is a discrete valuation on a field $F$, and
$X$ is a $d$-dimensional abelian variety over 
$F$ which has potentially good reduction at $v$. 
Let $G=G_{v,X}$. Suppose $\ell$ is a prime number not
equal to the residue characteristic of $v$. 
Then the action of $\Gal(F^s/F)$ on the $\ell$-adic Tate module
$V_\ell(X)$ induces an embedding
$$f : G \hookrightarrow \Sp_{2d}(\Q_\ell)$$
which satisfies the following properties.
\begin{enumerate}
\item[(i)] 
For every $\sigma \in G$, the coefficients of $P_{f(\sigma)}$
are integers which are independent of $\ell$.
If $X$ has an $F$-polarization of degree not divisible by $\ell$,
then one may choose $f$ so that its image lies in 
$\Sp_{2d}(\Z_\ell)$.
\item[(ii)] 
If either $(\ell,\#G)=1$ or $\ell > d+1$, then there exists
an embedding 
$$g : G \hookrightarrow \Sp_{2d}(\Z_\ell)$$
such that $P_{g(\sigma)} = P_{f(\sigma)}$ 
for every $\sigma \in G$.
\item[(iii)] 
If $\ell \ge 5$ then there exists an embedding 
$$h : G \hookrightarrow \Sp_{2d}(\F_\ell)$$
such that $P_{h(\sigma)} \equiv P_{f(\sigma)} \pmod{\ell}$ 
for every $\sigma \in G$.
\end{enumerate}
Further, if $\ell \ge 5$ then there exists an embedding 
$$G \hookrightarrow \Sp_{2d}(\Z_\ell)$$
(which does not necessarily ``preserve'' the characteristic
polynomials obtained from the embedding $f$).
\end{thm}

See \cite{Serre-Tate} for (i), and
see \cite{inertia} for the case $(\ell,\#G)=1$ of
(ii). The remainder of Theorem \ref{Gthm} follows from results
whose proofs will appear elsewhere (along with examples which
show that the results are sharp). Those results apply more
generally to measure how far an abelian variety (not necessarily
with potentially good reduction) is from having semistable reduction.
In some cases, these results apply to more general finite groups 
than those obtained as $G_{v,X}$'s.

\part{N\'eron models of abelian varieties with potentially good reduction}
\label{neronpart}

\section{Preliminaries}

In \cite{serrelem}, the following result was obtained as
a corollary of Theorem \ref{quasithm} above.

\begin{prop}[Theorem 6.10a of \cite{serrelem}]
\label{randm}
Suppose $\ell$ is a prime, $m$ and $r$ are positive integers, $\O$ is an
integral domain of characteristic zero with no non-zero infinitely 
$\ell$-divisible elements, $\ell\O$ is a maximal ideal of $\O$, $M$ is
a free $\O$-module of finite rank, and $A$ is an endomorphism of $M$ 
of finite multiplicative order such that
$(A - 1)^{m(\ell - 1)\ell^{r-1}} \in \ell^m\End(M)$.  
If $r > 1$, then the torsion
subgroup of $M/(A - 1)M$ is killed by $\ell^{r-1}$.
\end{prop}

\begin{prop}[see Proposition 6.1i and Corollary 7.1 of \cite{semistab}]
\label{pressred}
Suppose $X$ is a $d$-dimensional abelian variety over a field $F$, 
$v$ is a discrete valuation on $F$ with residue characteristic 
not equal to $2$, 
$\lambda$ is a polarization on $X$,   
$\widetilde{X}_2$ is a maximal isotropic subgroup of $X_2$ with respect to 
$e_{\lambda,2}$, $\lambda$ and the points of 
$\widetilde{X}_2$ are defined over an extension of $F$ which is 
unramified over $v$, 
${\bar v}$ is an extension of $v$ to a separable closure of $F$, 
and $\sigma \in \I({\bar v}/v)$.
Then $(\rho_{2,X}(\sigma) - I)^2 \in 2\M_{2d}(\Z_2)$, 
and $X$ has semistable 
reduction above $v$ over every totally ramified Galois 
(necessarily cyclic) extension of $F$ of degree $4$.
\end{prop}

Recall that $u$ denotes the unipotent rank of $X_v$, $a$ denotes the
abelian rank, and 
 $\Phi '$ denotes the prime-to-$p$ part of the group of 
connected components of the special
fiber of the N\'eron model of $X$ at $v$, where $p$ is the residue 
characteristic of the discrete valuation $v$. If $X$ has potentially
good reduction, then $\dim(X) = a + u$.

\begin{thm}[Theorem 7.5 of \cite{semistab}]
\label{neronmod}
Suppose $v$ is a discrete valuation on a field $F$ with strictly henselian 
valuation ring, $X$ is an abelian 
variety over $F$ which has potentially good reduction at $v$, and either 
\begin{enumerate}
\item[{(a)}] $n = 2$ and the points of $X_2$
are defined over $F$, or
\item[{(b)}] $n = 3$ or $4$, $\lambda$ is a polarization on $X$ 
defined over $F$, and 
the points of a maximal isotropic subgroup of $X_n$ with respect to 
$e_{\lambda,n}$ are defined over $F$.
\end{enumerate}
Suppose the residue characteristic $p$ ( $\ge 0$) of $v$ does
not divide $n$. 
Then $\Phi ' \cong (\Z/2\Z)^{2u}$
if $n = 2$ or $4$, and $\Phi ' \cong (\Z/3\Z)^u$ if $n = 3$.
\end{thm}

\begin{lem}
\label{neronlemma}
Suppose 
$v$ is a discrete valuation on a field $F$ such that 
the valuation ring is strictly henselian.
Suppose $X$ is an abelian 
variety over $F$ which has potentially good reduction at $v$,
and 
suppose $n$ is a positive integer not divisible by the residue
characteristic of $v$.
Let $\Phi_n$ denote the subgroup of 
$X_v/X_v^0$ of points of order dividing $n$.  
Then:
\begin{enumerate}
\item[(i)] $(X_v)_n \cong X_n(F)$,
\item[(ii)] $(X_v^0)_n \cong (\Z/n\Z)^{2a}$,  
\item[(iii)] $\Phi_n \cong (X_v)_n/(X_v^0)_n$,  and
\item[(iv)] if $X_n(F) \cong (\Z/n\Z)^{b}$,   
then $\Phi_n \cong (\Z/n\Z)^{b-2a}$.
\end{enumerate}
\end{lem}

\begin{proof}
By Lemma 2 of \cite{Serre-Tate}, the reduction map defines
an isomorphism of $X_n^\I$ onto $(X_v)_n$,
where $\I = \I({\bar v}/v)$ for some extension ${\bar v}$
of $v$ to $F^s$. 
Under our hypotheses on $v$, we have
$X_n^\I \cong X_n(F)$. Therefore, 
$(X_v)_n \cong X_n(F)$.
As shown in the proof of Lemma 1 of \cite{Serre-Tate},
$(X_v^0)_n \cong (\Z/n\Z)^{2a+t}$, where $t$ denotes
the toric rank of $X_v$.
Since $X$ has potentially good reduction at $v$, $t = 0$.
Since $X_v^0$ is $n$-divisible, we have 
$\Phi_n \cong (X_v)_n/(X_v^0)_n$. Part (iv) follows
easily from (i), (ii), and (iii).
\end{proof} 

\section{N\'eron models}
\label{main}

In Theorem \ref{neron} we generalize Theorem \ref{neronmod}
to the case of partial level $2$ structure. We can recover
Theorem \ref{neronmod}a as a special case. 
Recall that $u$ denotes the unipotent rank of $X_v$, 
$a$ denotes the abelian rank, and 
$\Phi '$ denotes the prime-to-$p$ part of the group of 
connected components of the special
fiber of the N\'eron model of $X$ at $v$, where $p$ is the residue 
characteristic of $v$ (with $\Phi '$ the full
group of components if $p = 0$). 

\begin{thm}
\label{neron}
Suppose $v$ is a discrete valuation on a field $F$, suppose
the valuation ring is strictly henselian, and suppose the residue field
has characteristic $p \ne 2$.
Suppose $(X, \lambda)$ is a $d$-dimensional polarized abelian 
variety over $F$, $X$ has potentially good reduction at $v$, and
the points of a maximal isotropic subgroup of $X_2$ with respect to 
$e_{\lambda,2}$ are defined over $F$.
Then:
\begin{enumerate}
\item[(i)]  $\Phi' \cong (\Z/2\Z)^{b-2a} = (\Z/2\Z)^{b+2u-2d}$, 
where $b$ is defined by $X_2(F) \cong (\Z/2\Z)^{b}$,
\item[(ii)]  $[X_2 : X_2(F)]\#\Phi' = 2^{2u}$, and
\item[(iii)]  $X$ has good reduction at $v$ if and only if
$\Phi' = \{0\}$ and $X_2 \subseteq X(F)$.
\end{enumerate}
\end{thm}

\begin{proof}
Let ${\bar v}$ be an extension of $v$ to a separable closure of $F$,
let $\I = \I({\bar v}/v)$, let $k$ be the residue field of $v$, and
let $\J$ be the first ramification group (i.e., $\J$ is trivial if
$p = 0$ and $\J$ is the pro-$p$-Sylow subgroup of $\I$ if $p > 0$). 
Suppose $q$ is a prime not equal to $p$, and let
$\Phi_q$ denote the $q$-part of the
group of connected components of the special
fiber of the N\'eron model of $X$. 
Since $X$ has potentially good reduction at $v$, $\rho_{q,X}(\sigma)$
has finite multiplicative order for every $\sigma \in \I$. 
Let $\tau$ be a lift to $\I$ of a generator of the pro-cyclic group
$\I/\J$. By \S11 of \cite{SGA} (see Lemma 2.1 of \cite{Lorenzini}), 
$$\Phi_q \text{ is isomorphic to the torsion subgroup of } 
T_q(X)^\J/(\tau - 1)T_q(X)^\J.$$
By Proposition \ref{pressred} and Remark \ref{ramifiedcyclic}, 
$X$ has semistable reduction (and therefore good reduction) above $v$ over 
a totally ramified cyclic Galois extension of $F$ of degree $4$. Therefore
$\I$ acts on $T_q(X)$ through a cyclic quotient of order $4$, so
$\rho_{q,X}(\sigma)^4 = I$ for every $\sigma \in \I$. Since $p \ne 2$,
we have $\rho_{q,X}(\sigma) = I$ for every $\sigma \in \J$. Therefore,
$T_q(X)^\J = T_q(X)$. If $q \ne 2$, then 
$T_q(X)/(\rho_{q,X}(\tau) - I)T_q(X)$
is torsion-free, so $\Phi_q$ is trivial. Further,
$$\Phi_2 \text{ is isomorphic to the torsion subgroup of } 
T_2(X)/(\tau - 1)T_2(X).$$ 
We have $(\rho_{2,X}(\tau) - I)^2 \in 2\M_{2d}(\Z_2)$, by Proposition \ref{pressred}.
By Proposition \ref{randm} with $\ell = 2$, $r = 2$, $m = 1$, and $\O = \Z_2$, 
$\Phi_2$ is annihilated by $2$.
Therefore, $\Phi'$ is an elementary abelian $2$-group.
By Lemma \ref{neronlemma},  $\Phi' \cong (\Z/2\Z)^{b-2a}$.
Parts (ii) and (iii) follow immediately.
Note that Theorem \ref{neronmod}a is a special case of Theorem
\ref{neron}.
\end{proof}

\begin{thm}
\label{neron3}
Suppose $v$ is a discrete valuation on a field $F$, suppose
the valuation ring is strictly henselian, and suppose the residue field
has characteristic $p \ne 3$.
Suppose $(X, \lambda)$ is a $d$-dimensional polarized abelian 
variety over $F$, $X$ has potentially good reduction at $v$, and 
the points of a maximal isotropic subgroup of $X_3$ with respect to 
$e_{\lambda,3}$ are defined over $F$. 
Then:
\begin{enumerate}
\item[(i)]  $X_3(F) \cong (\Z/3\Z)^{2d-u}$,  
\item[(ii)] $X$ has good reduction at $v$ if and only if
$X_3(F) = X_3$, and
\item[(iii)] $X$ has purely additive reduction at $v$ if and only if
$X_3(F) \cong (\Z/3\Z)^{d}$.
\end{enumerate}
\end{thm}

\begin{proof} 
By Theorem \ref{neronmod},
$\Phi' \cong (\Z/3\Z)^u$. Write $X_3(F) \cong (\Z/3\Z)^b$.
By Lemma \ref{neronlemma}, 
$\Phi' \cong (\Z/3\Z)^{b-2d+2u}$.
Therefore, $b = 2d - u$, and we obtain the desired result.
\end{proof}

\begin{thm}
\label{neron4}
Suppose $v$ is a discrete valuation on a field $F$ with strictly henselian 
valuation ring, $X$ is an abelian variety over $F$ which has potentially 
good reduction at $v$, the residue field has characteristic $p \ne 2$, 
and either 
\begin{enumerate}
\item[{(a)}] the points of $X_2$ are defined over $F$, or
\item[{(b)}] $\lambda$ is a polarization on $X$ defined over $F$, and 
the points of a maximal isotropic subgroup of $X_4$ with respect to 
$e_{\lambda,4}$ are defined over $F$.
\end{enumerate} 
Then
$$X_4(F) \cong (\Z/4\Z)^{2a} \times (\Z/2\Z)^{2u}.$$
In particular:
\begin{enumerate}
\item[(i)]  $X_2 \subseteq X_4(F) \subseteq X_4$, 
$[X_4 : X_4(F)] = 2^{2u}$, $[X_4(F) : X_2] = 2^{2a}$,
\item[(ii)] $X$ has good reduction at $v$ if and only if
$X_4(F) = X_4$, and
\item[(iii)] $X$ has purely additive reduction at $v$ if and only if
$X_4(F) = X_2$.
\end{enumerate}
\end{thm}

\begin{proof}
By Theorem \ref{neronmod},
we have $\Phi' \cong (\Z/2\Z)^{2u}$. By Lemma \ref{neronlemma}, 
we have a short exact sequence 
$$0 \to  (\Z/4\Z)^{2a} \to X_4(F) \to (\Z/2\Z)^{2u} \to 0.$$
Let $d = \dim(X)$. 
Since $X_4(F) \subseteq X_4 \cong (\Z/4\Z)^{2d}$,
we conclude that $X_4(F) \cong (\Z/4\Z)^{2a} \times (\Z/2\Z)^{2u}$. 
Note that $X_2 \cong (\Z/2\Z)^{2d} = (\Z/2\Z)^{2a + 2u}$. 
The rest of the result follows immediately.
\end{proof}

As an example, let $X$ be the elliptic curve defined by
the equation $y^2 = x^3 - 9x$, and let $F$ be the maximal unramified
extension of $\Q_3$.
Then $X_2(F) = X_2 = X_4(F)$, $X$ has additive and potentially
good reduction, and $\Phi' \cong (\Z/2\Z)^{2}$.

\begin{rems}
If $X$ has a polarization $\lambda$ of odd degree, then $X_2$ is 
a maximal isotropic subgroup of $X_4$ with respect to $e_{\lambda,4}$.

As stated in the Introduction, 
Theorems \ref{neron3}ii and \ref{neron4}ii are immediate
corollaries of Raynaud's criterion for semistable reduction. 

If $X$ has purely additive reduction, then $X_n(F) \cong
\Phi_n$ (see \cite{Lorenzini}).

Suppose $v$ is a discrete valuation on a field $F$,  
$X$ is an abelian variety over $F$ with potentially good reduction at $v$, 
the valuation ring is strictly henselian, $\ell = 2$ or $3$, and
the residue characteristic is not equal to $\ell$.
Then Theorem 6.1
of \cite{Edixhoven} implies that if $\Phi'$ is an elementary abelian
$\ell$-group,
then $\Phi'$ is a subgroup of $(\Z/2\Z)^{2u}$ if $\ell = 2$ or of 
$(\Z/3\Z)^{u}$ if $\ell = 3$.

For simplicity of exposition, we do not generalize the results of
\S \ref{main} (or the prerequisite results from \cite{semistab}, or
related results in \S 3 of \cite{connected})
to the setting of Theorem \ref{ssredconverse}, but leave such
generalizations as a straightforward exercise for the reader.
\end{rems}


\begin{thebibliography}{99}

\bibitem{BLR} S.\ Bosch, W.\ L\"utkebohmert, M.\ Raynaud, N\'eron models,
Springer, Berlin-Heidelberg-New York, 1990.
\bibitem{Edixhoven} B.\ Edixhoven, {\em On the prime-to-$p$ part of the
group of connected components of N\'eron models}, Comp.\ Math.\ {\bf 97} 
(1995), 29--49.
\bibitem{SGA} A.\ Grothendieck, {\em Mod\`eles de N\'eron et monodromie},
in Groupes de monodromie en g\'eometrie alg\'ebrique, SGA7 I,
A.\ Gro\-then\-dieck, ed., Lecture Notes in Math.\ {\bf 288}, Springer,
Berlin-Heidelberg-New York, 1972, pp.\ 313--523.
\bibitem{LenstraOort} H.\ W.\  Lenstra, Jr., F.\ Oort, {\em Abelian varieties
having purely additive reduction}, J.\ Pure and Applied Algebra {\bf 36}
(1985), 281--298.
\bibitem{Lorenzini} D.\ Lorenzini, {\em On the group of components of a
N\'eron model}, J.\ Reine Angew.\ Math.\ {\bf 445} (1993), 109--160.
\bibitem{Serre-Tate} J-P.\ Serre, J.\ Tate, {\em Good reduction of 
abelian varieties}, Ann.\ of Math.\ {\bf 88} (1968), 492--517. 
\bibitem{semistab} A.\ Silverberg, Yu.\ G.\ Zarhin, {\em Semistable reduction 
and torsion subgroups of abelian varieties}, 
Ann.\ Inst.\ Fourier {\bf 45}, no.~2 (1995), 403--420.
\bibitem{connected} A.\ Silverberg, Yu.\ G.\ Zarhin,
{\em Connectedness results for $\ell$-adic 
representations associated to abelian varieties}, 
Comp.\ math.\ {\bf 97} (1995), 273--284.
\bibitem{serrelem} A.\ Silverberg, Yu.\ G.\ Zarhin, {\em Variations on a 
theme of Minkowski and Serre}, 
J.\ Pure and Applied Algebra
{\bf 111}  (1996),  285--302.
\bibitem{degree} A.\ Silverberg, Yu.\ G.\ Zarhin, 
{\em Semistable reduction of abelian varieties over extensions of small degree},
J.\ Pure and Applied Algebra, to appear.
\bibitem{inertia} A.\ Silverberg, Yu.\ G.\ Zarhin, {\em Subgroups of inertia
groups arising from abelian varieties}, J. Algebra, to appear.
\bibitem{etale} A.\ Silverberg, Yu.\ G.\ Zarhin, 
{\em \'Etale cohomology and reduction of abelian varieties}, preprint.
\bibitem{Silverman} J.\ H.\ Silverman, {\em The N\'eron fiber of
abelian varieties with potential good reduction}, Math.\ Ann.\ {\bf 264}
(1983), 1--3.
\bibitem {WeilAV} A.\ Weil, Vari\'et\'es ab\'eliennes et courbes 
alg\'ebriques, Hermann, Paris, 1948.

\end{thebibliography}
\end{document}